\documentclass[12pt]{article}
\usepackage{amssymb}
\usepackage{amsmath}
\usepackage{amsfonts}
\usepackage{tabularx}
\usepackage[francais,english]{babel}
\usepackage[latin1]{inputenc}
\usepackage{epsfig}
\usepackage{theorem}
{\theoremstyle{break}
{\theorembodyfont{\normalfont}
\newtheorem{rem}{Remarque }}}

\begin{document}
%
\begin{center}
{\Large    {\sc         Classification supervisée en grande dimension.
\newline
Application à l'agrément de conduite automobile    } }

\bigskip
Jean-Michel Poggi \& Christine Tuleau

\medskip
{\it
  Laboratoire de Mathématique\ --\ U.M.R. C
8628, ``Probabilités, Statistique et Modélisation'',
Université Paris-Sud, Bât.~425, 91405 Orsay Cedex, France

(et Univ. Paris 5)}

jean-michel.poggi@math.u-psud.fr christine.tuleau@math.u-psud.fr
\end{center}
\bigskip

{\small {\bf \noindent Résumé}

\noindent Ce travail est motivé par un problème réel
appelé l'objectivation. Il consiste à expliquer
l'agrément de conduite au moyen de critères ``physiques'',
issus de signaux mesurés lors d'essais. Nous suggérons une
approche pour le problème de la sélection des variables
discriminantes en tentant de tirer profit du caractère
fonctionnel des données. Le problème est mal posé, au
sens où le nombre de variables explicatives est très
supérieur à la taille de l'échantillon. La démarche
procède en trois étapes : un prétraitement des signaux
incluant débruitage par ondelettes, recalage et
synchronisation, une réduction de la taille des signaux par
compression dans une base d'ondelettes commune, et enfin
l'extraction des variables utiles au moyen d'une stratégie
incluant des applications successives de la méthode CART.\\

{\it \noindent Mots clés :} CART, Classification, Discrimination,
Ondelettes}\\

{\small {\bf  \noindent Abstract}

\selectlanguage{francais}
\noindent This work is motivated by a real world problem:
objectivization. It consists of explaining the subjective drivability using physical criteria coming from signals measured
during experiments. We suggest an approach for the discriminant
variables selection trying to take advantage of the functional
nature of the data. The problem is ill-posed, since the number of
explanatory variables is hugely greater than the sample size. The
strategy proceeds in three steps: a signal preprocessing,
including wavelet denoising and synchronization, dimensionality
reduction by compression using a common wavelet basis, and finally
the selection of useful variables using a stepwise strategy
involving successive applications of the CART method.\\
\selectlanguage{english}

}

\section{Introduction}
\selectlanguage{francais}
Ce travail est motivé par un problème réel appelé
l'objectivation. Il consiste à expliquer l'agrément de
conduite traduisant un confort ressenti relativement à une
prestation donnée, par exemple le comportement de la boîte de
vitesses lors de la phase de mise en mouvement d'un véhicule, au
moyen de critères ``physiques'', c'est-à-dire de variables
issues de signaux (comme une vitesse, des couples ou encore la position
de pédales) mesurés lors d'essais. Il s'agit d'utiliser cette quantification 
pour en tenir compte lors de la phase de conception du véhicule. Il s'inscrit dans la
continuité de travaux menés par Renault portant sur la
prestation décollage à plat pour un groupe moto-propulseur
à boîte de vitesses robotisée (cf. Ansaldi \cite{Ansaldi}). \\
\selectlanguage{english}

\selectlanguage{francais}
Dans cet article, nous  développons une approche alternative pour
le problème de la  sélection des variables discriminantes en
tentant de plus tirer profit du caractère fonctionnel des
données. De ce point de vue, ce travail peut être
rapproché de l'analyse des données fonctionnelles. Citons
Deville \cite{Deville}, Dauxois et Pousse \cite{Dauxois} pour les
travaux pionniers dans les années 70. Plus récemment, on peut citer par exemple,
Leurgans \textit{et al.} \cite{Leurgans}, Hastie \textit{et al.}
\cite{Hastie2} et ces dernières années, Ferraty, Vieu \cite{Ferraty}, Ferré \textit{et al.} (\cite{Ferre}, 
\cite{Villa}), Rossi et Conan-Guez \cite{Rossi}, Biau \textit{et al.} \cite{Biau} ainsi que le texte de
synthèse de Besse, Cardot \cite{Besse}. En outre les deux livres de
Ramsay, Silverman \cite{Ramsay}, \cite{Ramsay2} constituent une
ressource précieuse.\\
Dans ce travail, nous préférons utiliser la méthode CART particulièrement adaptée pour la sélection de variables.\\
\selectlanguage{english}

\selectlanguage{francais}
Comme cela est classique dans de nombreuses applications où les
variables explicatives sont des courbes, le problème industriel
qui nous occupe est mal posé, au sens où le nombre de
variables explicatives est très supérieur à la taille de
l'échantillon. L'un des exemples typiques de telles situations est
fourni par les données  d'expression du génome. On trouvera
dans Dudoit {\it et al.} \cite{Dudoit} la présentation de ce
problème et de diverses méthodes de classification
supervisée actuellement en compétition. On pourra aussi
consulter Vannucci {\it et al.} \cite{Vannucci} pour la situation où
les variables explicatives sont des spectres, ce qui est classique en
chimiométrie.\\
\selectlanguage{english}

\selectlanguage{francais}
Structurellement le problème industriel qui nous intéresse présente
une particularité supplémentaire : nous disposons non pas d'une seule
variable explicative qui est une courbe mais d'un grand nombre de
variables fonctionnelles parmi lesquelles il faut choisir les plus
influentes. Notre approche s'intéresse donc à un double problème de
sélection : celle des variables fonctionnelles d'une part, et d'autre
part pour chacune de ces courbes, la sélection de bons descripteurs
discriminants.\\
\selectlanguage{english}

\selectlanguage{francais}
La démarche adoptée procède en trois étapes et utilise deux
outils fondamentaux que sont d'une part la méthode des ondelettes
(cf. Misiti {\it et al.} \cite{Misiti}) et d'autre part la méthode de
classification non linéaire CART (cf. \cite{Breiman}). Les trois
étapes sont constituées d'un prétraitement des signaux
(incluant débruitage par ondelettes, recalage et synchronisation),
d'une réduction de la dimension par compression dans une base
d'ondelettes commune puis de l'extraction et sélection des variables
utiles au moyen d'une stratégie incluant des applications
successives de la méthode CART.\\
\selectlanguage{english}

\selectlanguage{francais}
Le plan de l'article est le suivant. Après cette introduction, le
paragraphe 2 présente le contexte de l'application : le
problème et les données. Dans le paragraphe 3, la démarche
adoptée est détaillée. Enfin le paragraphe 4 regroupe quelques
éléments de conclusion.
\selectlanguage{english}

\section{Le contexte applicatif}

\subsection{Le problème}
\selectlanguage{francais}
\noindent
La campagne d'essais réalisée par Renault (cf. Ansaldi
\cite{Ansaldi}) a conduit à faire varier les facteurs suivants
: le réglage de la boîte de vitesses, les conditions de
roulage et les pilotes. Lors de ces essais, ont été
mesurés d'une part l'agrément du pilote et d'autre part
des données objectives consistant dans le relevé, à
l'aide de capteurs, de plusieurs signaux temporels.\\
\selectlanguage{english}

\selectlanguage{francais}
\noindent
Précisons quelques éléments de terminologie utiles dans la suite. On
appelle ``produit'' un élément de

$\{produits\}\ =\ \{conditions\ de\ roulage\}\ \times\ \{3\
r\acute{e}glages\ de\ la\ boite\ de\ vitesses\}$ \\
\noindent
où \\
$\{conditions\ de\ roulage\}\ =\{2\ charges\}\ \times\ \{2\ angles\
p\acute{e}dale\} \times\ \{2\ vitesses\ p\acute{e}dale\}$
\selectlanguage{english}

\selectlanguage{francais}
\noindent
ce qui conduit au plus à 24 produits (12 pour chacune des charges).\\
\selectlanguage{english}

\selectlanguage{francais}
\noindent
On appelle ``essai'' un élément de $\{essais\}\ =\ \{7\
pilotes\}\ \times\ \{24\ produits\}$
\selectlanguage{english}

\selectlanguage{francais}
\noindent
conduisant à un maximum de 168 essais. \\
\selectlanguage{english}

\selectlanguage{francais}
\noindent
Les essais à 140 kg de charge ont été menés
séparément des essais à 280 kg de charge. Pour chaque charge, 
6 produits parmi les 12 possibles ont été testés :
4 pilotes ont comparé par paires ces 6 produits. Après analyse des
résultats, 114 essais à 140 kg et 118 essais à 280 kg ont
été retenus. Chacun de ces essais est représenté par un
ensemble de 21 variables fonctionnelles qui correspondent aux signaux
mesurés par les capteurs durant l'expérience.\\
\selectlanguage{english}

\selectlanguage{francais}
\noindent
L'étude menée dans \cite{Ansaldi} s'articule autour de trois phases :
\selectlanguage{english}
\begin{itemize}
\item l'association d'un agrément à chacun des produits.\\
Pour chaque paire d'essais, le pilote précisait
son essai préféré. A partir de ces données
de comparaisons par paires et à l'aide d'une
méthode inspirée du ``multidimensional scaling'' (voir la thèse de Favre \cite{Favre})
sont obtenus un classement des
produits par pilote et un agrément consensuel à toute la
population des pilotes, par charge. Cet agrément associe à
un produit un rang de satisfaction (le rang 1 étant celui du
produit le plus apprécié);

\item l'extraction de critères puis sélection par analyse
  discriminante. \\
A partir des signaux mesurés, de très nombreux critères
  sont générés puis, au moyen d'une analyse discriminante
  linéaire arborescente dite par moindres écarts 
(c'est-à-dire basée sur un critère $L^1$), un petit nombre
  d'entre eux expliquant l'agrément, sont extraits;

\item le calcul d'intervalles de tolérance. \\
Pour chacun des critères pertinents, un intervalle qui
maximise l'agrément sous certaines contraintes sur les
produits, est construit (ce point constitue d'ailleurs la
contribution majeure de la thèse d'Ansaldi \cite{Ansaldi}).
\end{itemize}

\selectlanguage{francais}
\noindent
On se concentre, dans cet article, sur la deuxième étape en
utilisant une approche plus fonctionnelle. Bien sûr, on ne considère
que les données issues de la phase 1 qui sont seules détaillées dans
le paragraphe suivant. L'agrément est le rang consensuel attribué 
à chacun des 6 produits testés. Ceci conduit à 
un problème de discrimination, au lieu d'un problème de 
régression avec une variable à expliquer ordinale discrète.\\ 
\selectlanguage{english}

\selectlanguage{francais}
\noindent
Dans la suite, ne seront considérés que les essais à 140 kg de charge 
(pour les essais à 280 kg de charge, la démarche est identique et les 
résultats obtenus dans l'étude \cite{Ansaldi} sont semblables). \\
\selectlanguage{english}

\subsection{Les données}
\selectlanguage{francais}
\noindent
Les données sont constituées des couples $ ((X^{j}_{i})_{1
  \leq{j} \leq{J}}, Y_{i})_{1 \leq{i} \leq{n} }$, où $n=114$ et
$J=21$, et :
\selectlanguage{english}
\begin{enumerate}
\item [-] $Y_{i}$ représente le rang attribué au produit
  testé lors de l'essai $i$;
\item [-] $X^{j}_{i}$ représente la $j^{\grave{e}me}$ variable
  fonctionnelle mesurée lors de l'essai $i$ et est le signal $\{X^{j}_{i}(t)\}_{t
    \in{T_{i}}}$ où $T_{i}$ est la grille temporelle
  régulière propre à l'essai $i$.
\end{enumerate}

\selectlanguage{francais}
\noindent
Autrement dit, pour chacun des essais, on dispose de l'agrément et
de 21 signaux (on parlera dans la suite, suivant le contexte, de
signaux comme de variables fonctionnelles ou encore de courbes) pour
la plupart d'environ 1000 points (en fait ils comportent entre 600 et
5000 points). Ces variables fonctionnelles sont principalement des positions, 
des vitesses, des accélérations, des couples et des régimes moteur, cependant pour des raisons de confidentialité la nature des variables ne peut pas être indiquée de façon plus précise. Notons que la fréquence d'échantillonnage de 250 Hz est 
la même pour tous les essais et correspond à une haute résolution temporelle.\\
\selectlanguage{english}

\selectlanguage{francais}
\noindent
La distribution de l'agrément $Y$, après regroupement en $5$ modalités, est donnée par les fréquences $33 \%, 17 \%, 17 \%, 18 \%, 15 \%  $. Seulement $5$ modalités, et non $6$, sont prises en considération, deux produits ayant obtenu le même agrément.\\
\selectlanguage{english}

\begin{figure}[!h]
\begin{center}
\begin{tabularx}{13cm}{X}
\multicolumn{1}{c}{
\includegraphics[width=11cm,height=6.5cm]{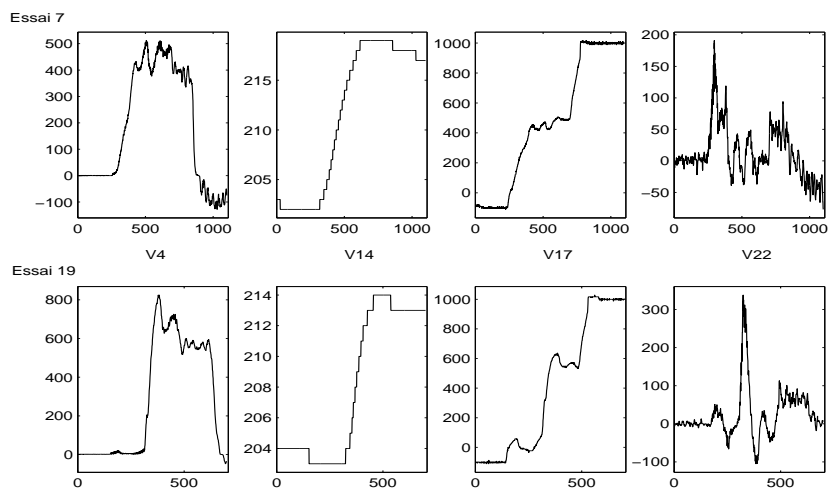}}\\
\caption{\footnotesize Pour les essais 7 et 19, les quatre variables
  fonctionnelles $X^{j}$ correspondant à $j=4, 14, 17, 22$, notées
  simplement V4, V14, V17 et V22. Elles sont observées  sur une grille
  temporelle propre à l'essai et présentent des  caractéristiques
  temporelles variées.} \label{sigbruts}
\end{tabularx}
\end{center}
\end{figure}

\selectlanguage{francais}
\noindent On trouve dans la Figure \ref{sigbruts}, les quatre
variables fonctionnelles $X^{j}$ correspondant à $j=4$, $14, 17,
22$ pour les essais $7$ et $19$.\\

\selectlanguage{francais}
\selectlanguage{english}

\selectlanguage{francais}
\noindent L'examen des graphiques permet de formuler quelques
remarques préliminaires concernant ces variables fonctionnelles :
\selectlanguage{english}
\begin{itemize}
\item \selectlanguage{francais} elles sont observées sur une grille temporelle propre à l'essai,
  ce qui nécessitera des recalages temporels;
\selectlanguage{english}
\item \selectlanguage{francais} elles peuvent être d'allure générale et d'ordre de grandeur très
  différents, à la fois pour un même essai mais aussi au travers des
  différents essais, ce qui impliquera des recalages en ordonnée des
  courbes;
\selectlanguage{english}
\item \selectlanguage{francais} elles présentent des caractéristiques temporelles très
  différentes, par exemple le rapport signal sur bruit, élevé en
  général, peut s'avérer modéré comme dans le cas de la variable 22 ou
  encore l'être localement comme c'est le cas pour ces quatre
  variables sauf la variable 14 qui est une fonction constante par
  morceaux. Il est clair qu'un débruitage, sans être en général
  crucial, peut s'avérer utile;
\selectlanguage{english}
\item \selectlanguage{francais} la forme générale est souvent simple et peu de paramètres ou peu
  d'événements semblent suffisants pour la caractériser. Ceci
   permet d'espérer à la fois une caractérisation économe
  des variables fonctionnelles ainsi qu'une compression efficace.
\selectlanguage{english}
\end{itemize}

\begin{rem}
\selectlanguage{francais}
\noindent
La variabilité, entre les essais, des durées d'observation et celle des 
amplitudes des signaux mesurés, résultent des différences de conditions 
de roulage et de l'exécution plus ou moins scrupuleuse des consignes par les pilotes. 
\selectlanguage{english}
\end{rem}


\section{La démarche}
\selectlanguage{francais}
\noindent
Le cadre général dans lequel on se place est celui de la sélection de
variables dans un problème de discrimination, et consiste à
construire une fonction, génériquement notée $F$ dans la suite, pour
prédire $Y$ à l'aide de :
$$\hat Y = F(X^{1},...,X^{J})$$
Dans cette perspective, il sera utile de sélectionner
parcimonieusement les variables fonctionnelles qui peuvent expliquer
l'agrément, puis pour chacune d'elles, de ne retenir qu'un
très faible nombre d'aspects la décrivant, pour des raisons
évidentes de robustesse.\\
\selectlanguage{english}

\selectlanguage{francais}
\noindent Autrement dit, on cherche à sélectionner ce que nous
appelons dans ce contexte, des critères notés $C^{j_k}$, déduits
des $X^{j}$, de façon à prédire convenablement $Y$ par
:

$$\hat Y = F(C^{j_1},...,C^{j_K})$$ avec $K<<J$, typiquement de
l'ordre de 5 pour l'application industrielle.\\

\selectlanguage{english}

\selectlanguage{francais}
\noindent 
Rappelons que dans le cadre de l'objectivation, il ne
s'agit pas d'expliquer au mieux l'agrément en utilisant toutes les
informations disponibles, comme par exemple les conditions de roulage, 
qui ont un impact certain, mais de l'expliquer partiellement en se 
restreignant exclusivement à des variables
déduites des signaux mesurés de fa\c con à pouvoir remonter à des paramètres de conception du véhicule. \\

\selectlanguage{english}

\selectlanguage{francais}
\noindent La démarche adoptée procède en trois étapes :
\selectlanguage{english}
\begin{itemize}
\item \selectlanguage{francais} un prétraitement des signaux, incluant débruitage par
  ondelettes, recalage et synchronisation;
\selectlanguage{english}
\item \selectlanguage{francais} une réduction de la taille des signaux par compression dans
  une base d'ondelettes commune;
\selectlanguage{english}
\item \selectlanguage{francais} l'extraction des variables utiles au moyen d'une stratégie
  pas à pas procédant par des applications successives de la
  méthode CART.
\selectlanguage{english}
\end{itemize}

\selectlanguage{francais}
\noindent Détaillons successivement chacune de ces trois phases.
\selectlanguage{english}

\subsection{Prétraitement des signaux}
\selectlanguage{francais}
\noindent Les données $X^{j}_{i}\ =\ \{X^{j}_{i}(t)\}_{t
\in{T_{i}}}$ sont prétraitées de façon d'une part,
à les débruiter individuellement c'est-à-dire pour un essai et
une variable fonctionnelle donnés et, d'autre part,  à les
rendre plus homogènes au moyen de recalages.
\selectlanguage{english}

\selectlanguage{francais}
\subsubsection{Tronquer les signaux}
\noindent Avant ces deux traitements, on isole une phase qui est
la seule à être directement déduite de connaissances externes
propres au problème. En effet, en dépit de consignes clairement 
définies, les durées d'enregistrement et les dates des différentes 
étapes de l'essai ne sont pas synchrones. Néanmoins, on peut définir
deux événements à réaligner : le ``vrai'' début de 
l'essai et sa ``vraie'' fin qui sont lisibles au travers des variables 
fonctionnelles 8 et 21. Ces deux événements correspondent physiquement 
au démarrage réel du véhicule et à la définition de la fin de l'essai.\\ 
\selectlanguage{english}

\selectlanguage{francais}
\begin{figure}[!h]
\begin{center}
\begin{tabularx}{14cm}{X}
\multicolumn{1}{c}{
\includegraphics[width=11cm,height=7cm]{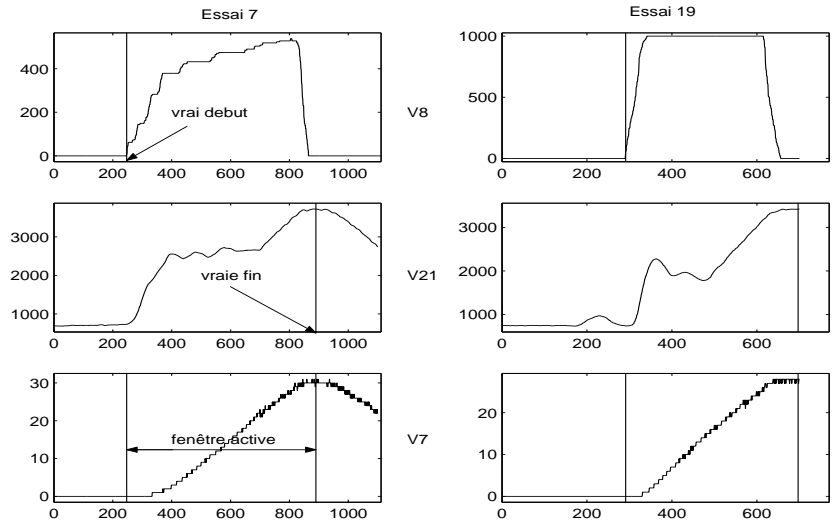}}\\
\caption{\footnotesize Pour les essais 7 et 19, les trois variables
  fonctionnelles $X^{j}$ correspondant à $j=8, 21, 7$ notées
  simplement V8, V21 et V7 sur le graphique. Les
  deux premières servent de marqueur au ``vrai'' début de l'essai
  et sa ``vraie'' fin, respectivement. La période utile de l'essai est
  visualisée sur les graphes de la variable fonctionnelle 7 par la
  portion de signal située entre les deux instants matérialisés par
  des lignes verticales.} \label{fenactive}
\end{tabularx}
\end{center}
\end{figure}
\selectlanguage{english}

\selectlanguage{francais}
\noindent 
On trouve dans la Figure \ref{fenactive}, trois variables fonctionnelles
$X^{j}$ correspondant à $j=8, 21, 7$ pour les 
essais $7$ et $19$. Les deux premières servent de marqueur du ``vrai''
début de l'essai et de sa ``vraie'' fin, respectivement. La
période utile de l'essai est visualisée sur les graphes de la
variable fonctionnelle 7 par la portion de signal située entre les
deux instants matérialisés par des lignes verticales. Bien sûr,
ces instants varient en fonction de l'essai.

\selectlanguage{francais}
\noindent
Pour l'essai $i$, on note $\widetilde{T}_{i}$ la grille ${T_{i}}$
convenablement tronquée aux extrémités.
\selectlanguage{english}

\subsubsection{Débruiter les signaux}
\selectlanguage{francais}
\noindent A $i$ et $j$ fixés, le signal mesuré est
contaminé par un bruit de capteur. Bien sûr, il convient de
l'éliminer avant tout traitement de type recalage ou interpolation
des données, qui conduirait à les modifier et donc altérer la
nature stochastique du bruit qui affecte le signal utile. Comme
l'atteste la Figure \ref{sigbruts}, la régularité locale de
celui-ci peut beaucoup varier au cours du temps, il convient donc
d'utiliser des techniques de débruitage adaptatives en espace.
C'est le cas de celles basées sur les méthodes d'ondelettes (cf.
Donoho, Johnstone \cite{Donoho} pour l'un des articles fondateurs,
Vidakovic \cite{Vidakovic} pour un large tour d'horizon de ces
méthodes et Misiti {\it et al.} \cite{Misiti} pour une
introduction aisée).
\\
\selectlanguage{english}

\selectlanguage{francais}
\noindent On considère le modèle suivant, usuel en traitement
statistique du signal et réaliste dans cette application :
$$\forall t \in{\widetilde{T}_{i}},\ \ \ X^{j}_{i}(t)\ =\
f^{j}_{i}(t)\ +\ \eta^{j}_{i}(t)$$ où $\{\eta^{j}_{i}(t)\}_{t \in
  \widetilde{T}_{i}}$ est un bruit blanc. Dans ce cadre, le débruitage
consiste à décomposer le signal dans une base d'ondelettes, à seuiller
les coefficients de détail de façon à éliminer essentiellement ceux
attribuables au bruit puis à reconstruire un signal débruité
constitué de la somme d'une approximation lisse et de détails à
diverses échelles correspondant aux fluctuations rapides du signal
utile.
\selectlanguage{english}

\begin{figure}[!h]
\begin{center}
\begin{tabularx}{13cm}{X}
\multicolumn{1}{c}{
\includegraphics[width=11cm,height=8cm]{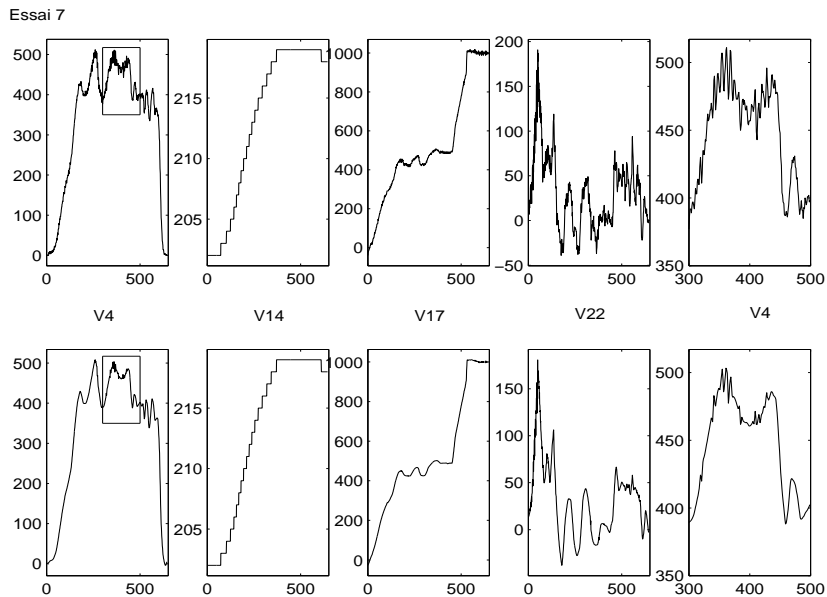}}\\
\caption{\footnotesize Pour l'essai 7, en haut de la figure les
quatre signaux $X^{j}_{7}(t)$ \hspace{2cm} ($j=4, 14, 17, 22$) et, en bas, leurs
versions débruitées. Dans les deux derniers graphiques à droite,
un zoom sur une portion du premier signal permet d'apprécier la
qualité du débruitage par ondelettes, à la fois efficace pour
débruiter les parties lisses tout en préservant les composantes à
haute fréquence du signal utile.} \label{debruitage}
\end{tabularx}
\end{center}
\end{figure}

\selectlanguage{francais}
\noindent On obtient ainsi une estimation $ \{\hat{f}^{j}_{i}(t)\}_{t
  \in{\widetilde{T}_{i}}}$, ou encore un signal débruité
$\{\hat{X}^{j}_{i}(t)\}_{t \in{\widetilde{T}_{i}}}$.
\noindent La Figure \ref{debruitage} présente les résultats
obtenus après débruitage par ondelettes des quatre variables
fonctionnelles montrées en Figure \ref{sigbruts}. La méthode 
utilise l'ondelette de Daubechies presque symétrique d'ordre 
$4$, un niveau de décomposition entre $3$ et $5$ 
(suivant les signaux) et le seuillage dit "universel" 
(cf. Donoho et Johnstone \cite{Donoho}).\\
\selectlanguage{english}

\selectlanguage{francais}
\selectlanguage{english}

\selectlanguage{francais}
\noindent Comme on peut le remarquer, le débruitage par ondelettes
permet de supprimer de façon satisfaisante le bruit tout en préservant
les composantes à haute fréquence du signal utile.
\selectlanguage{english}

\subsubsection{Synchroniser et normaliser les signaux}
\selectlanguage{francais}
\noindent
L'objectif de cette étape est d'éliminer la dépendance en
$i$ de la grille temporelle. On procède pour chaque signal, tout
d'abord à un recalage linéaire en temps en ramenant la grille
$\widetilde{T}_{i}$ sur l'intervalle $[0,1]$. Puis, on effectue une
interpolation linéaire du signal, suivie d'un échantillonnage
pour se ramener à la grille régulière à $m$ points de
$[0,1]$ (ici on fixe $m=512$, valeur largement suffisante pour des durées de 
fenêtres actives comprises entre $300$ et $700$ observations). Un instant 
dans cette nouvelle ``unité'' de
temps s'interprète comme la proportion de la durée de l'essai
écoulée.
\selectlanguage{english}

\selectlanguage{francais}
\noindent
On dispose donc de $\{\widetilde{X}^{j}_{i}(t)\}_{t \in{T}}$, sur la
grille fixe $T=\{\frac{1}{m},...,\frac{m-1}{m},1 \}$.\\
\selectlanguage{english}

\selectlanguage{francais}
\noindent Enfin, pour éliminer certains effets d'échelle, en
partie liés aux conditions de roulage, les signaux sont
normalisés en ordonnée.\\
\selectlanguage{english}

\begin{rem}
\selectlanguage{francais}
\par \noindent Un autre prétraitement consiste à effectuer un
recalage non linéaire en alignant pour tout $j$, les $n$
signaux à l'aide de marqueurs convenablement choisis (cf.
Bigot \cite{Bigot}). Ceci amènerait à considérer le
problème plus sous un aspect de classification de formes.
Cependant, cela serait extrêmement lourd et engendrerait une
difficulté quant à la remontée dans le temps d'origine
en particularisant de nouveau les variables fonctionnelles, et
limiterait l'interprétation.\\ 
\selectlanguage{english}

\selectlanguage{francais}
\noindent 
En revanche, cela permettrait de poursuivre un objectif plus
ambitieux consistant à augmenter l'homogénéité à $Y$
fixé, en mettant au point le recalage pour chaque modalité de la
réponse. \\
\selectlanguage{english}

\selectlanguage{francais}
\noindent
Mentionnons que des méthodes de recalage temporel intermédiaires 
entre la solution adoptée et celle-ci sont envisageables, comme par 
exemple le type de méthode de recalage décrit dans \cite{Ramsay} 
qui cherche à rapprocher des fonctions de leur moyenne.\\
\selectlanguage{english}
\end{rem}

\begin{rem}
\selectlanguage{francais}
\noindent De manière implicite, dans la suite du travail (mais
aussi dans les travaux antérieurs menés dans ce contexte par
Renault), les essais sont considérés comme des
réplications indépendantes. Des classifications non
supervisées et des ACP fonctionnelles (cf. Ramsay, Silverman
\cite{Ramsay}) permettent de corroborer raisonnablement l'idée que
les effets dus au pilote et aux conditions de roulage sont
négligeables devant les autres facteurs de variabilité.
\selectlanguage{english}
\end{rem}

\subsection{Compression des signaux}
\selectlanguage{francais}
\noindent
A l'issue de la phase de prétraitement, on dispose donc pour chaque
essai, de $J=21$ signaux débruités, de $m=512$ points. Chacun
de ces signaux peut donc être représenté dans une base
d'ondelettes ou de paquets d'ondelettes par très peu de
coefficients (cf. Mallat \cite{Mallat} et Coifman, Wickerhauser
\cite{Coifman2}). Il suffit, par exemple, pour un signal donné, de
sélectionner les coefficients les plus grands en valeur absolue,
exploitant ainsi la capacité des ondelettes à concentrer l'énergie
d'un signal (pour des classes très larges de signaux), en un très
petit nombre de ses grands coefficients d'ondelettes. \\ Le
problème est ici de choisir, variable fonctionnelle par variable 
fonctionnelle, une base commune à tous les essais
pour les représenter de façon compacte. Pour déterminer
une base commune de décomposition, on peut se restreindre à un
petit nombre de bases différentes comme les espaces d'approximation en
ondelettes de résolution de plus en plus grossière. Comme $512=2^9$,
seule une demie douzaine de bases, l'ondelette étant choisie 
(ici on utilise l'ondelette de Daubechies presque symétrique d'ordre $4$), 
sont à mettre en compétition. Le choix peut être :
\selectlanguage{english}
\begin{itemize}
\item \selectlanguage{francais} effectué indépendamment de la variable $Y$ et guidé par la
  définition d'un critère de qualité comme par exemple la
  moyenne de l'erreur d'approximation du signal par sa projection
  convenablement pénalisé.
\selectlanguage{english}

\selectlanguage{francais}
\noindent
Afin de déterminer le niveau de décomposition de chacun des signaux $j$, 
on considère le critère $EQ_j(p)$ lié à l'énergie et défini 
comme suit :
\selectlanguage{english}
\begin{itemize}
\item \selectlanguage{francais} pour une variable fonctionnelle $j$ et pour un individu $i$, soit 
$X^j_i(t)$ le signal d'origine et $A^j_{i,p}(t)$ le signal reconstruit 
à partir des coefficients d'approximation du niveau $p$;
\selectlanguage{english}
\item \selectlanguage{francais} on définit l'erreur de la variable fonctionnelle $j$ par 
$$EQ_j(p)=\sum_{i=1}^{114} \Vert X^j_i(t)-A^j_{i,p}(t) {\Vert}^2$$
\selectlanguage{english}
\end{itemize}

\begin{rem}
\selectlanguage{francais}
\noindent Notons que, lorsque le niveau de décomposition $p$ augmente, le nombre de coefficients et la qualité d'approximation diminuent. Le choix du niveau de décomposition résulte d'un compromis entre le nombre de coefficients retenus et la qualité d'approximation.
\selectlanguage{english}
\end{rem}

\selectlanguage{francais}
\noindent
Le choix du niveau de décomposition de la variable $j$ consiste alors 
à déterminer la plus petite valeur de $p$ pour laquelle on détecte 
un changement de pente ``suffisant'' dans le graphe de 
$(p,EQ_j(p))_{1 \leq p \leq 9}$ et à \^{o}ter $1$ à titre conservatoire.\\
\selectlanguage{english}

\selectlanguage{francais}
\begin{figure}[!h]
\begin{center}
\begin{tabularx}{13cm}{X}
\multicolumn{1}{c}{
\includegraphics[width=11cm,height=7.5cm]{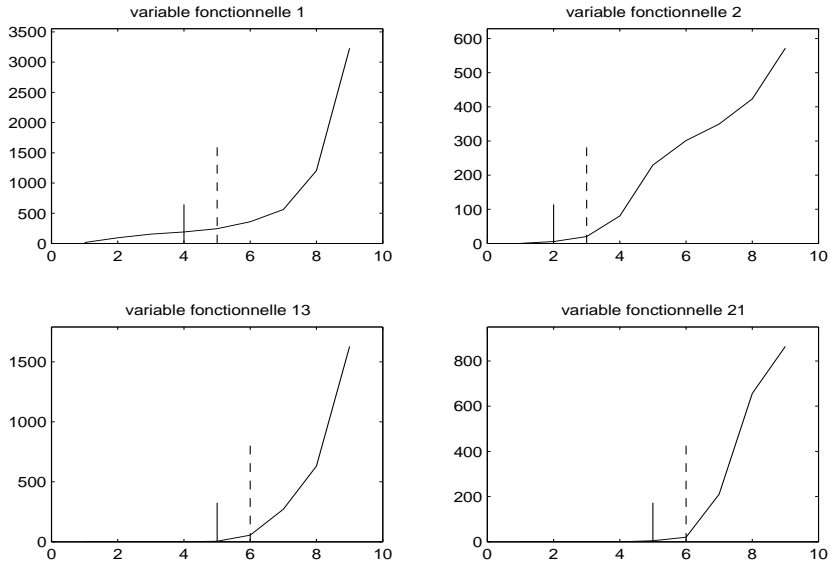}}\\
\caption{\footnotesize Pour les variables fonctionnelles $1,2,13$ et $21$, on 
représente $(p,EQ(p))_{1 \leq p \leq 9}$ en trait plein, la plus petite valeur 
de $p$ pour laquelle on détecte un changement de pente ``suffisant'' en traits 
pointillés et cette valeur ôtée de $1$ en traits pleins.} \label{energie}
\end{tabularx}
\end{center}
\end{figure}
\selectlanguage{english}

\selectlanguage{francais}
\noindent
La Figure \ref{energie} esquisse la fa\c con dont le niveau de 
décomposition lors de la compression par ondelettes est déterminé 
pour chacun des signaux.\\
\selectlanguage{english}

\item \selectlanguage{francais} basé sur un critère dépendant de la variable $Y$, comme
  par exemple l'erreur de classification d'un arbre CART (voir
  paragraphe suivant).
\selectlanguage{english}
\end{itemize}
\vspace*{0.5cm}
\selectlanguage{francais}
\noindent
L'emploi d'une procédure inspirée du premier choix ci-dessus avec
recherche d'une cassure dans la répartition moyenne de l'énergie,
conduit à retenir majoritairement 16 coefficients et donc à
réduire $\mathbb{R}^{J\ \times\ m}$ à
$\mathbb{R}^{\sum{m_{j}}}$ avec $\sum{m_{j}}\ \thickapprox\ 300$ ou
$400$ suivant la stratégie adoptée pour comprimer une variable
fonctionnelle (d'ailleurs non discriminante) dont les fluctuations à
haute fréquence sont significatives.\\
\selectlanguage{english}

\begin{figure}[!h]
\begin{center}
\begin{tabularx}{13cm}{X}
\multicolumn{1}{c}{
\includegraphics[width=11cm,height=8cm]{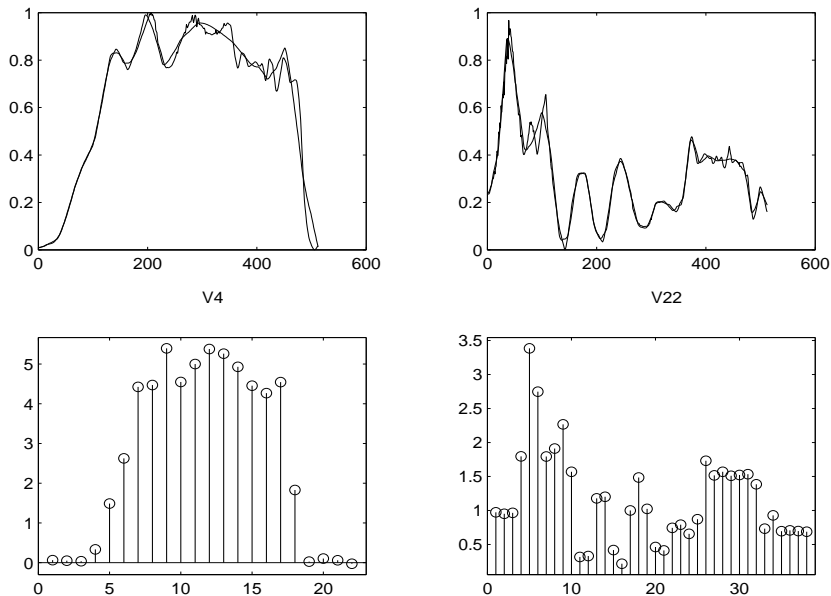}}\\
\caption{\footnotesize Pour l'essai 7 et pour les deux variables
  correspondant à $j=4, 22$ : en haut, le signal après compression
  superposé au signal original (prétraité), en bas les coefficients
  d'approximation associés.} \label{compression}
\end{tabularx}
\end{center}
\end{figure}

\selectlanguage{francais}
\noindent La Figure \ref{compression} présente pour deux variables
fonctionnelles, les résultats obtenus après compression par
ondelettes : le signal après compression superposé au signal
prétraité est représenté dans le premier graphique, le second (en
dessous) contient les coefficients d'approximation associés.
Ceux-ci peuvent, bien sûr, être de taille différente puisque le
niveau de décomposition retenu dépend de la variable considérée.\\
\selectlanguage{english}

\selectlanguage{francais}
\selectlanguage{english}

\selectlanguage{francais}
\noindent Les deux graphiques du haut de la Figure \ref{compression}
contiennent, pour l'essai 7 et pour deux variables fonctionnelles
différentes, le signal après compression superposé au signal
d'origine. Ils sont très proches bien que représentés par peu
de coefficients. En effet, les deux graphiques du bas de la figure
contiennent les coefficients d'approximation associés aux
représentations comprimées. Ainsi, la forme des graphiques du haut et du
bas de la figure se ressemblent sauf aux extrémités de
l'axe des abscisses à cause d'extra-coefficients, engendrés par les
prolongements appliqués aux signaux dans les calculs des
coefficients par la transformée en ondelettes discrète (voir
\cite{Misiti}).
\selectlanguage{english}

\begin{rem}
\selectlanguage{francais}
\noindent
Signalons que la connexion entre les développements sur des bases
orthogonales d'ondelettes de processus stochastiques et les
décompositions issues de la transformée discrète en ondelettes est
donnée, par exemple, dans Amato {\it et al.} \cite{am03}.
\selectlanguage{english}
\end{rem}

\begin{rem}
\selectlanguage{francais}
\noindent Une autre approche associant plus étroitement les phases
de compression et de sélection des variables discriminantes est
proposée par Coifman, Saito \cite{Coifman}. Il s'agit de choisir une
base optimale, parmi les bases associées à une
décomposition en paquets d'ondelettes, en maximisant la
séparation entre classes.
\selectlanguage{english}
\end{rem}

\selectlanguage{francais}
\noindent Elle n'est pas retenue ici, une voie médiane est
empruntée : des gains massifs en compression sont obtenus même au
prix d'un politique de sélection un peu conservative de façon à ne
pas trop obérer la phase suivante qui fera le choix des variables
les plus discriminantes. On note $C^{j}=(C^{j,1},...,C^{j,K_{j}})$ le paquet des $K_{j}$
coefficients associés à la variable
fonctionnelle $X^{j}$.
\selectlanguage{english}

\subsection{Sélection de variables par CART}
\selectlanguage{francais}
\noindent
A la fin de l'étape précédente, il y a une
réduction de la dimension de l'espace des variables, mais elle
demeure insuffisante puisque l'on dispose de $114$ individus à
comparer à $300$ ou $400$ variables.
\selectlanguage{english}

\selectlanguage{francais}
\noindent
Les nouvelles données ainsi construites sont donc :
$ (((C^{j,k}_{i})_{1 \leq{k} \leq{K_{j}}})_{1 \leq{j} \leq{J}},
Y_{i})_{1 \leq{i} \leq{n} }$.\\

\noindent On propose une procédure pas à pas basée sur
la méthode CART. Celle-ci permet d'ajuster aux données, un
modèle additif du type $Y=F((C^{j,k})_{j,k})$ où $F$ est
additive et plus précisément constante sur des polyèdres dont les
côtés sont parallèles aux axes, sous la forme d'un arbre dyadique
de décision. On peut se reporter au livre de Breiman {\it et al.}
\cite{Breiman} les fondateurs de la méthode ou
Hastie {\it et al.} \cite{Hastie} pour un rapide aperçu. Dans la
suite, on considère l'erreur de classification définie comme
usuellement mais en pénalisant les fausses classifications par le
truchement de la matrice de coût définie par $\Gamma(k,k')=|k-k'|$,
définition qui découle naturellement du fait que $Y$ est une variable
ordinale discrète.\\
\selectlanguage{english}

\begin{figure}[!h]
\begin{center}
\begin{tabularx}{13cm}{X}
\multicolumn{1}{c}{
\includegraphics[width=10cm,height=5cm]{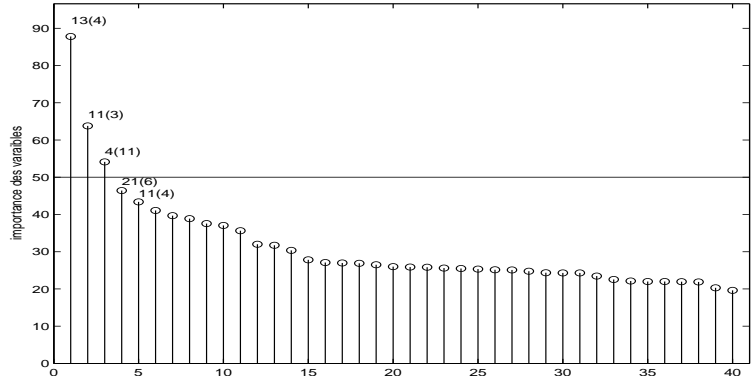}}\\
\caption{\footnotesize Importance des variables calculée sur
le modèle $M^{j_0}$ sélectionné précédemment et
sélection finale des trois critères dont les importances
ressortent nettement en tête.} \label{impof}
\end{tabularx}
\end{center}
\end{figure}

\selectlanguage{francais}
\noindent La procédure est présentée ci-dessous en cinq phases :
\selectlanguage{english}
\begin{enumerate}
\item \selectlanguage{francais}
Pour chaque $j$, on construit l'arbre CART $A^{j}$ expliquant $Y$
par le paquet de coefficients $C^{j}$ et on sélectionne, au
moyen de l'importance des variables au sens de Breiman {\it et
al.} \cite{Breiman} (voir aussi \cite{Ghattas} et \cite{Ghattas2}), 
le paquet des coefficients utiles, noté
$\tilde{C}^{j}$, en seuillant l'importance comme illustré dans
la Figure \ref{coefut}.
\selectlanguage{english}

\selectlanguage{francais}
\begin{figure}[!h]
\begin{center}
\begin{tabularx}{13cm}{X}
\multicolumn{1}{c}{
\includegraphics[width=11cm,height=8cm]{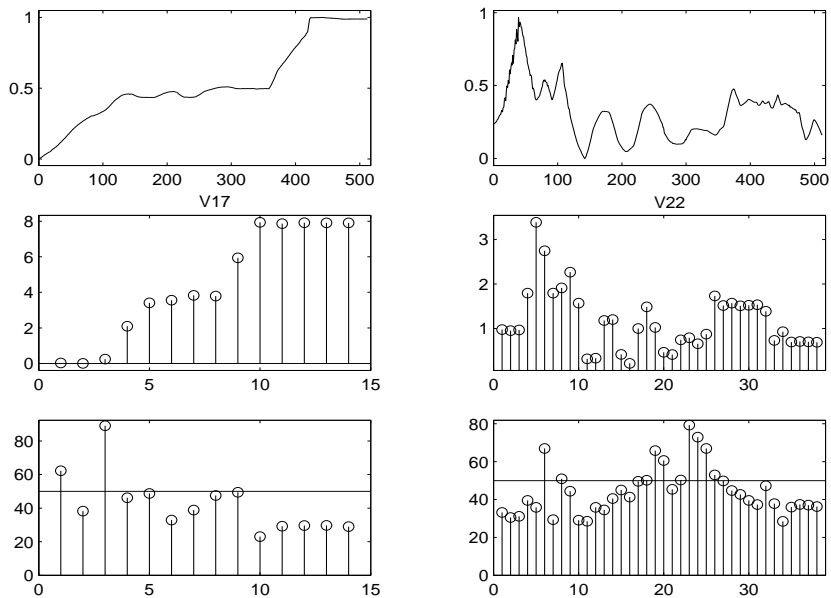}}\\
\caption{\footnotesize Pour les essais 7 et 19, et pour les
variables correspondant à \hspace{2cm} $j=17,22 $, en haut les signaux
prétraités, au milieu le paquet $C^{j}$ des coefficients
d'approximation de niveau retenu  et en bas l'importance de chacun
de ces coefficients. Les coefficients utiles constituant
$\tilde{C}^{j}$ sont ceux dont l'importance dépasse le seuil.}
\label{coefut}
\end{tabularx}
\end{center}
\end{figure}
\selectlanguage{english}

\selectlanguage{francais}
\noindent
On peut noter que les pics dans les graphes de l'importance des
variables correspondent non pas seulement, à des marqueurs
significatifs de la forme du signal mais bien à des
événements significatifs discriminants.
\selectlanguage{english}
\item \selectlanguage{francais}
On en déduit un ordre sur les ``nouvelles'' variables
fonctionnelles (c'est-à-dire sur les paquets $(\tilde{C}^{j})_j$)
au moyen de l'erreur de classification, évaluée par validation
croisée, commise par l'arbre $A^{j}$ (voir Figure \ref{coutva}).
\selectlanguage{english}

\selectlanguage{francais}
\begin{figure}[!h]
\begin{center}
\begin{tabularx}{13cm}{X}
\multicolumn{1}{c}{
\includegraphics[width=8cm,height=5cm]{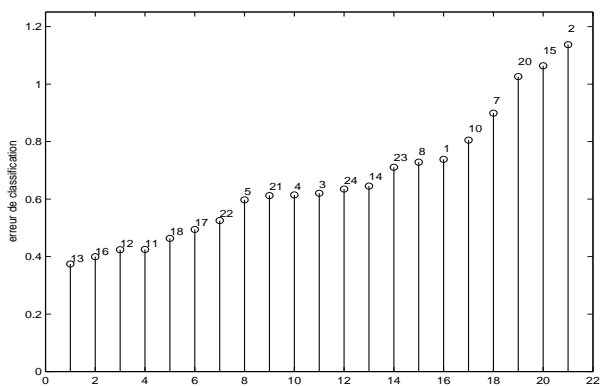}}\\
\caption{\footnotesize L'erreur de classification évaluée par
validation croisée des arbres
$A^{j}$, de la meilleure à la pire. Elle fluctue dans un rapport
de 1 à 3. Cet ordre sur les ``nouvelles'' variables fonctionnelles
est celui qui sera utilisé pour les invoquer pas à pas.}
\label{coutva}
\end{tabularx}
\end{center}
\end{figure}
\selectlanguage{english}

\item \selectlanguage{francais}
On construit une suite ascendante $(M^{j})_{j}$ d'au plus $J=21$ modèles CART
emboîtés, en invoquant et en testant les paquets de variables
$\tilde{C}^{j}$, pas à pas, suivant l'ordre précédemment
obtenu. Autrement dit, $M^{j}$ explique $Y$ par l'ensemble de paquets
de coefficients $(\tilde{C}^{l})_{{l} \leq{j}}$ privés des paquets qui se sont révélés, après test, comme insuffisamment informatifs.
\selectlanguage{english}

\selectlanguage{francais}
\begin{figure}[!h]
\begin{center}
\begin{tabularx}{13cm}{X}\multicolumn{1}{c}{
\includegraphics[width=10cm,height=5cm]{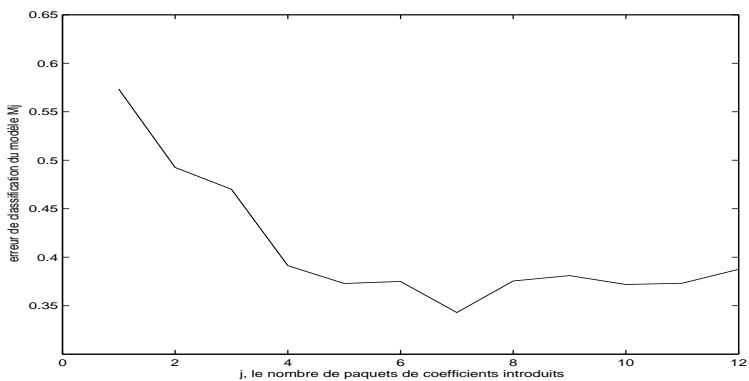}}\\
\caption{\footnotesize L'erreur de classification du modèle
$M^{j}$ évaluée par validation
croisée, en fonction de $j$, le nombre de paquets de
coefficients introduits.} \label{coutmod}
\end{tabularx}
\end{center}
\end{figure}
\selectlanguage{english}

\item \selectlanguage{francais}
On sélectionne ensuite les variables fonctionnelles
pertinentes en choisissant celles définissant le modèle
$M^{j_0}$ minimisant l'erreur de classification. L'allure de
celle-ci (cf. Figure \ref{coutmod}) est attendue : elle
décroît d'abord fortement avant de lentement
croître lorsque les variables introduites n'apportent plus rien à
la discrimination.
\selectlanguage{english}

\item \selectlanguage{francais}
Enfin, en calculant l'importance des variables explicatives du
modèle $M^{j_0}$ : les coefficients $\{\tilde{C}^j, j \in M^{j_0}\}$ 
et en retenant la tête de ce classement, on sélectionne
les critères pertinents (voir Figure \ref{impof}).
\selectlanguage{english}
\end{enumerate}

\selectlanguage{francais}
\selectlanguage{english}

\selectlanguage{francais}
\noindent
Une première façon de procéder, très dépendante du problème, consiste à ne retenir que les $5$ premières variables, $5$ étant le nombre souhaité de critères. On obtient alors un arbre dont l'erreur de validation croisée est de $24$ sur $114$ pour $12$ erreurs apparentes (c'est-à-dire l'erreur de resubstitution).\\
Une alternative consiste à considérer l'erreur de validation croisée sur la suite de modèles emboîtés induite par l'ordre issu du calcul de l'importance des variables. On sélectionne alors le modèle dont l'erreur est la plus faible.\\
\selectlanguage{english}

\selectlanguage{francais}
\begin{table}[!h]
\begin{center}
\begin{tabularx}{12cm}{X}
\multicolumn{1}{c}{
\includegraphics[width=11cm,height=2cm]{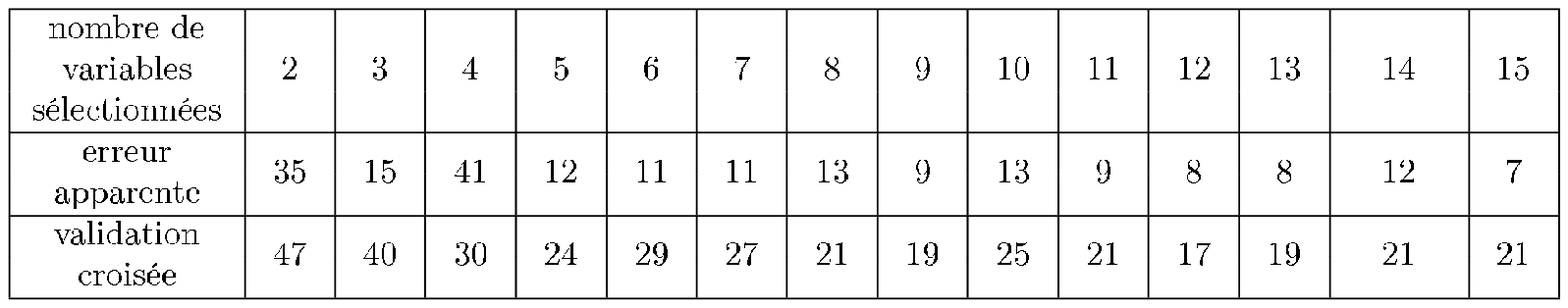}}\\
\caption{\label{resulf} Nombre d'erreurs commises sur l'échantillon d'apprentissage en fonction du nombre de variables retenues.}
\end{tabularx}
\end{center}
\end{table}
\selectlanguage{english}

\begin{figure}[!h]
\begin{center}
\begin{tabularx}{13cm}{X}
\multicolumn{1}{c}{
\includegraphics[width=7cm,height=12cm,angle=270]{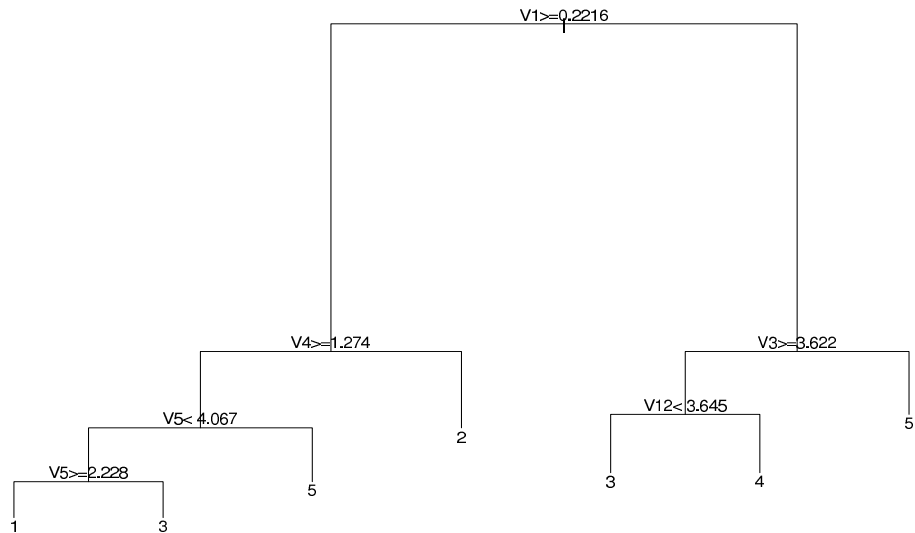}}\\
\caption{\footnotesize Arbre dont l'erreur de classification, évaluée par validation croisée, est la plus faible.} \label{arbf}
\end{tabularx}
\end{center}
\end{figure}

\selectlanguage{francais}
\noindent
La Table \ref{resulf} donne, pour les modèles de cette suite dont le nombre de variables est inférieur à $15$, l'erreur apparente et l'erreur de validation croisée. Le meilleur modèle est celui comportant $12$ variables. L'erreur commise est de $17$ sur $114$ ($15\%$) et l'erreur apparente de $8$ sur $114$ ($7\%$), ce qui est très satisfaisant.\\
Enfin, si l'on examine l'arbre CART construit en se restreignant à ces $12$ variables (cf. Figure \ref{arbf}), il est intéressant de noter que $5$ variables seulement étiquettent les n\oe uds de l'arbre et $4$ d'entre elles sont en tête du classement fourni par la Figure \ref{impof}.
\selectlanguage{english}

\selectlanguage{francais}
\selectlanguage{english}

\begin{rem}
\selectlanguage{francais}
\noindent
Terminons par une remarque générale dont la portée méthodologique est cruciale.\\ 
Un inconvénient classique de l'usage des arbres de
classification est leur instabilité, c'est-à-dire que le
classifieur construit peut fluctuer ``beaucoup'' pour des
``petites'' variations de l'échantillon d'apprentissage (cf.
Hastie {\it et al.} \cite{Hastie}). Un remède désormais classique
à cette propriété indésirable est d'utiliser le bagging qui permet
de stabiliser la prédiction en utilisant non pas un classifieur
mais l'agrégé d'un ensemble de classifieurs construits par
rééchantillonnage bootstrap de l'échantillon d'apprentissage (voir
Ghattas \cite{Ghattas2}).\\
\selectlanguage{english}

\selectlanguage{francais}
\noindent Suivant cette idée (voir Ghattas \cite{Ghattas}),
l'importance des variables et l'erreur de classification sont
évaluées par rééchantillonnage. Plus précisément,
pour la phase 1, on considère la moyenne des importances des
variables calculées sur des arbres obtenus par rééchantillonnage 
$n$ pour $n$, des $114$ observations. Pour
l'estimation de l'erreur de classification, elle est évaluée par
validation croisée grâce à un schéma de découpage en 10 de
l'échantillon puis stabilisée en randomisant cette phase de
découpage.
\\
\selectlanguage{english}
\end{rem}

\section{Conclusion}
\selectlanguage{francais}
\noindent
Du point de vue de l'application, les critères qui
ressortent comme les plus discriminants sont associés à quatre variables fonctionnelles. Parmi eux, deux sont
très proches des critères obtenus par la méthode basée sur la méthode
discriminante linéaire et deux sont nouveaux et considérés par les experts comme intéressants. Il faut noter que
dans notre cas, ces critères ont été obtenus sans intégrer de
connaissances {\textit{a priori}}, sauf dans la phase de troncature
de la grille temporelle des observations. Signalons cependant que
les conditions d'arrêt dépendent de seuils fixés pour
le moment en fonction de l'application.\\
\selectlanguage{english}

\selectlanguage{francais}
\noindent Complémentairement à ce travail, des avancées concernent l'étude
théorique de pénalités adéquates pour faire de la sélection de
variables dans des contextes voisins (cf. Sauvé, Tuleau
\cite{Sauve}). Typiquement il s'agit d'utiliser une approche par
sélection de modèle ``à la Birgé-Massart''
(cf. Barron, Birgé, Massart \cite{BirgeMassart}) pour sélectionner des variables dans un
modèle de régression non linéaire, au moyen d'applications
répétées de la méthode CART. Des résultats de type
inégalités oracles permettent de préciser la forme des
pénalités convenables et peuvent suggérer des alternatives au
choix ad-hoc effectués ici.\\
\selectlanguage{english}

\noindent {\large{\bf Remerciements}}\\
\selectlanguage{francais}
\noindent Les auteurs remercient la Direction de la Recherche de
Renault d'avoir mis à leur disposition les données
relatives aux essais qui motivent ce travail et, en particulier,
Nadine Ansaldi pour les discussions associées. Cette collaboration se poursuit actuellement dans le cadre d'un
contrat de recherche entre le laboratoire de mathématiques d'Orsay
et la Direction de la Recherche de Renault. \\
\selectlanguage{english}

\selectlanguage{francais}
\noindent
En outre, les auteurs remercient les deux rapporteurs anonymes de leurs remarques et suggestions qui ont contribué grandement à clarifier et améliorer la première version du manuscrit.\\
\selectlanguage{english}
\selectlanguage{francais}
\bibliographystyle{plain}
\bibliography{bibli}

\end{document}